
\input amstex.tex
\documentstyle{amsppt}
\magnification=\magstep1
\hsize=12.5cm
\vsize=18cm
\hoffset=1cm
\voffset=2cm
\def\DJ{\leavevmode\setbox0=\hbox{D}\kern0pt\rlap
{\kern.04em\raise.188\ht0\hbox{-}}D}
\footline={\hss{\vbox to 2cm{\vfil\hbox{\rm\folio}}}\hss}
\nopagenumbers 
\font\ff=cmr8
\def\txt#1{{\textstyle{#1}}}
\baselineskip=13pt
\def\hf{{\textstyle{1\over2}}}
\def\a{\alpha}\def\b{\beta}
\def\d{{\,\roman d}}
\def\e{\varepsilon}

\def\g{\gamma} \def\G{\Gamma}

\def\s{\sigma}

\def\={\;=\;}

\def\zt{\zeta(\hf+it)}

\def\no{\noindent}  
\def\R{\Re{\roman e}\,} \def\I{\Im{\roman m}\,} \def\s{\sigma}
\def\z{\zeta} 
\def\no{\noindent} 
\def\e{\varepsilon}
 
\def\no{\noindent} 
\def\e{\varepsilon}

\def\no{\noindent} 
\font\teneufm=eufm10
\font\seveneufm=eufm7
\font\fiveeufm=eufm5
\newfam\eufmfam
\textfont\eufmfam=\teneufm
\scriptfont\eufmfam=\seveneufm
\scriptscriptfont\eufmfam=\fiveeufm
\def\mathfrak#1{{\fam\eufmfam\relax#1}}

\font\tenmsb=msbm10
\font\sevenmsb=msbm7
\font\fivemsb=msbm5
\newfam\msbfam
\textfont\msbfam=\tenmsb
\scriptfont\msbfam=\sevenmsb
\scriptscriptfont\msbfam=\fivemsb
\def\Bbb#1{{\fam\msbfam #1}}

\def \NN {\Bbb N}
\def \CC {\Bbb C}

\def\rightheadline{{\hfil{\ff
Sums over ordinates of zeta-zeros}\hfil\tenrm\folio}}

\def\leftheadline{{\tenrm\folio\hfil{\ff
Aleksandar Ivi\'c }\hfil}}
\def\emptyheadline{\hfil}
\headline{\ifnum\pageno=1 \emptyheadline\else
\ifodd\pageno \rightheadline \else \leftheadline\fi\fi}

\topmatter
\title ON CERTAIN SUMS OVER ORDINATES OF ZETA-ZEROS \endtitle
\author   Aleksandar  Ivi\'c \endauthor
\address{ \bigskip
Aleksandar Ivi\'c, Katedra Matematike RGF-a
Universiteta u Beogradu, \DJ u\v sina 7, 11000 Beograd,
Serbia (Yugoslavia). \bigskip}
\endaddress
\dedicatory
Bulletin CXXI de l'Acad\'e\-mie Serbe des Sciences et des 
Arts - 2001, Classe des Sciences math\'ematiques et naturelles, 
Sciences math\'ematiques No. {\bf26}, pp. 39-52.
\enddedicatory
\email {\tt aleks\@ivic.matf.bg.ac.yu, 
aivic\@rgf.rgf.bg.ac.yu} \endemail
\abstract
{Let $\g$ denote imaginary parts of complex zeros
of $\z(s)$. Certain sums over the $\g$'s are evaluated,
by using the function $G(s) = \sum_{\g > 0}\g^{-s}$
and other techniques. Some integrals involving the function
$S(T)$ are also considered.}
\endabstract
\endtopmatter

\vglue 1cm
\heading{\bf1. The function $G(s)$}
\endheading
\noindent
Define, for $\s = \R s > 1$,
$$
G(s) = \sum_{\g > 0}\g^{-s},\eqno(1.1)
$$
where $\g$ denotes ordinates of complex zeros of the Riemann zeta-function
$\z(s)$. The aim of this note is to provide the (unconditional)
study of $G(s)$ and some applications to the evaluation
 of sums over the $\g$'s
and some related integrals. The function $G(s)$ is mentioned,
in a perfunctory way, in the work of Chakravarty [2] and in more detail by
Delsarte [5]. A related zeta-function, namely
$$
\sum_{\g > 0}\g^{-s}\sin(\a\g)\qquad(\a > 0),
$$
was studied by Fujii [6], but its properties are different from
the properties of $G(s)$, and we shall not consider it here. Both
Chakravarty and Delsarte (as well as Fujii) assume the Riemann Hypothesis
(that all complex zeros of $\z(s)$ satisfy $\R s = \hf$, RH for short) 
in dealing with $G(s)$. Delsarte [5] obtains its analytic continuation 
to $\CC$ under the RH. This will be obtained later in Section 3 by an 
argument which is different from Delsarte's, who employed a sort of 
a modular relation to deal with $G(s)$.

\medskip
To begin the study of $G(s)$ we need some notation. 
As usual, let the function
$$
N(T) \= \sum_{0<\g\le T}1
$$
count the number of positive imaginary parts of all complex zeros
which do not exceed $T$. We have (see [4, Chapter 15] or [13, Section 9.3])
$$\eqalign{
N(T) &= \sum_{0<\g\le T}1 = {1\over\pi}\vartheta(T) + 1 + S(T),\cr
\vartheta(T) &= \I\left\{\log\Gamma({\txt{1\over4}}+\hf iT)\right\}
- \hf T\log\pi,
\cr}
$$
where $\vartheta(T)$ is continuously differentiable, 
and if $T$ is not an ordinate of a zero
$$ 
S(T) = {1\over\pi}\arg\z(\hf+iT) = {1\over\pi}\I\left\{
\log\z(\hf+iT)\right\} \ll \log T. \eqno(1.2)
$$
Here the argument of $\z(\hf+iT)$ is
obtained by continuous variation along the straight lines joining the
points 2, $2+iT$, $\hf+iT$, starting with the value 0. If $T$
is  an ordinate of a zero, then $S(T) = S(T+0)$.\
 
It is clear then that the series in (1.1) converges absolutely for
$\s > 1$, and to  obtain its analytic continuation to the
region $\s \le 1$ we use Stirling's formula
for the gamma-function (see [8]) and write the formula for $N(T)$ as
$$
N(T) = {T\over2\pi}\log{T\over2\pi} - {T\over2\pi} + {7\over8}
+ S(T) + f(T), \; f(T) \ll {1\over T},\; f'(T) \ll {1\over T^2}.\eqno(1.3)
$$
Since the smallest positive ordinate of a zeta-zero is $14.13\ldots$, we have
$$
\eqalign{
G(s) &= \int_1^\infty x^{-s}\d N(x) =
\int_1^\infty x^{-s}\left\{{1\over2\pi}\log({x\over2\pi})\d x + 
\d\left(S(x) + f(x)\right)\right\}
\cr&
= {1\over2\pi}\left({x^{1-s}\over1-s} \log({x\over2\pi})\Big|_1^\infty
- \int_1^\infty {x^{1-s}\over1-s}\cdot{\d x\over x}\right)
\cr&
+ x^{-s}\left(S(x) + f(x)\right)\Big|_1^\infty 
+ s\int_1^\infty\left(S(x) + f(x)\right)x^{-s-1}\d x.\cr}
$$
In view of the bounds in (1.2) and (1.3) the last integral is seen
to converge absolutely.
Thus by the principle of analytic continuation we have, for $\s > 0$,
$$
G(s) = {1\over2\pi(s-1)^2} - {\log2\pi\over2\pi(s-1)} + C_1 
+ s\int_1^\infty\left(S(x) + f(x)\right)x^{-s-1}\d x,\eqno(1.4)
$$
where $C_1$ is a suitable constant. A relation similar to (1.4)
was established by Chakravarty [3, p. 490]. Further analytic continuation
will follow by integrating by parts the last integral.
This will give, for $\s > -1$,
$$\eqalign{
G(s) &= {1\over2\pi(s-1)^2} - {\log2\pi\over2\pi(s-1)} + C_1 
\cr&+ s\int_1^\infty f(x)x^{-s-1}\d x + s(s+1)\int_1^\infty
\int_1^x S(u)\d u\cdot x^{-s-2}\d x,\cr}\eqno(1.5)
$$
since we have the bound (see [13])
$$
\int_0^T S(t)\d t \;=\;O(\log T).\eqno(1.6)
$$
It follows that (1.5) gives
$$
G(s) \ll t^2\qquad(\s > -1,\,|t| \ge t_0).\eqno(1.7)
$$
Hence by convexity (the Phragm\'en-Lindel\"of principle, see [8]) we have
$$
G(s) \ll_\e |t|^\e(1+|t|^{1-\s})\qquad(\s > -1,\,|t| \ge t_0),\eqno(1.8)
$$
since $G(s) \ll 1$ for $\s > 1$. A sharper bound than (1.8),
at least for $0 \le \s\le 1$,  can be obtained as follows. We have
(initially    for $\s > 1$, then by analytic continuation for $\s > 0$)
$$
G(s) = \sum_{0<\g\le X}\g^{-s} + \sum_{\g> X}\g^{-s}
= \sum\nolimits_1(s,X) + \sum\nolimits_2(s,X),\eqno(1.9)
$$
say. The function $\sum_1(s,X)$ is entire, and we have by partial summation
(since $N(T) \ll T\log T$)
$$
\sum\nolimits_1(s,X) \ll X^{1-\s}\log X +\log^2X\qquad(0 \le \s \le 1).
$$
Henceforth we suppose that $T \le t \le 2T$ and we shall choose
$X = X(T) \,(\ge 2)$ appropriately a little later. Integration by
parts gives
$$
\eqalign{
\sum\nolimits_2(s,X) &= \int_X^\infty x^{-s}\left({1\over2\pi}
\log\left({x\over2\pi}\right)\d x + \d(S(x)+f(x))\right)\cr&
= {X^{1-\s}\over s-1}\cdot{1\over2\pi}
\log\left({X\over2\pi}\right) + {1\over s-1}\int_X^\infty {x^{-s}
\over2\pi}\d x \cr&
+ O(X^{-\s}\log^2X) + s\int_X^\infty (S(x)+f(x))
x^{-s-1}\d x.\cr}\eqno(1.10)
$$
This gives
$$
\eqalign{
&\;G(s)\ll \cr&\ll X^{1-\s}\log X + X^{-\s}\log^2X + X^{1-\s}|t|^{-2}\log X +
X^{-\s}|t|\log\log X + \log^2X
\cr&
\ll |t|^{1-\s}\log|t| + \log^2|t|\quad(X = T,\,0<\s<1).\cr}
$$
Therefore by continuity we obtain a sharpening of (1.8)
for $0 \le\s\le1$, namely
$$
G(s) \ll |t|^{1-\s}\log|t| + \log^2|t|\qquad(0 \le \s \le 1,\,|t| \ge
t_0 > 0).\eqno(1.11)
$$
In estimating the last integral in (1.10) we used the Cauchy-Schwarz
inequality for integrals and the mean square bound for $S(t)$
(see (4.3)).

\bigskip
\heading{\bf2. Mean square estimates for $G(s)$}
\endheading
\noindent
We pass now to mean square estimates for $G(s)$, for which
as usual we expect to smoothen the irregularites of the integrand.
If $0\le\s\le1$, then we can write
$$
\eqalign{
\int_T^{2T}|G(\s+it)|^2\d t &\ll \int_{T}^{2T}
|\sum\nolimits_1(\s+it,X)|^2\d t \,+ \cr&
+\,\int_{T}^{2T}
|\sum\nolimits_2(\s+it,X)|^2\d t = I_1(T) + I_2(T),\cr}\eqno(2.1)
$$
say, where $\sum_1$ and $\sum_2$ are defined by (1.9). 
To bound $I_1(T)$ we use the mean value theorem for Dirichlet
polynomials (see e.g., [8, Th. 5.2]) in the form
$$
\int_0^T\bigl|\sum_{n\le N}a_nn^{-it}\bigr|^2\d t =
T\sum_{n\le N}|a_n|^2 + O\bigl(\sum_{n\le N}n|a_n|^2\bigr).\eqno(2.2)
$$
If $0  < \g_1 \le \g_2 \le \cdots\;$ denote positive ordinates of zeta
zeros, then we can write
$$
\sum\nolimits_1(\s+it,X) = \sum_{\g_n\le X}\g_n^{-\s}\g_n^{-it},
\qquad \g_n \asymp n\log n.
$$
Hence with $X = T$ and $a_n = \g_n^{-\s}$ we obtain from (2.2)
$$
I_1(T) \ll \Bigg\{\aligned T\qquad
&(\s > \hf),\\ T\log^2T\qquad&(\s = \hf),\\
T^{2-2\s}\log T \qquad&(\s < \hf).\endaligned\eqno(2.3)
$$
To bound $I_2(T)$, we recall Parseval's formula for Mellin transforms
(see [12]) in the form
$$
\int_0^\infty f(x)g(x)x^{2\s-1}\d x = {1\over2\pi i}\int_{(\s)}
F(s)\overline{G(s)}\d s,\eqno(2.4)
$$
provided that
$$
H(s) = \int_0^\infty h(x)x^{s-1}\d x,\quad x^{\s-{1\over2}}h(x) \in
L^2(0,\infty)
$$
with $h(x) = f(x)$ or $h(x) = g(x)$. As usual $\int_{(c)}$ 
denotes $\lim_{T\to\infty}\int_{c-iT}^{c+iT}$. From (2.4) one obtains
$$
\int_1^\infty f(x)g(x)x^{1-2\s}\d x = {1\over2\pi i}\int_{(\s)}
F^*(s)\overline{G^*(s)}\d s,\eqno(2.5)
$$
provided that
$$
H^*(s) = \int_1^\infty h(x)x^{-s}\d x,\quad x^{{1\over2}-\s}h(x) \in
L^2(0,\infty)
$$
with $h(x) = f(x)$ or $h(x) = g(x)$. Setting in (2.5) $f(x) = g(x)$
if $a \le x \le b \, (1 \le a < b)$ and $f(x) = 0$ otherwise, it follows that
$$
\int_T^{2T}\Bigl|\int_a^b g(x)x^{-\s-it}\d x\Bigr|^2\d t
\le 2\pi\int_a^b g^2(x)x^{1-2\s}\d x.\eqno(2.6)
$$
Applying (1.10), (2.6) and (4.3) we obtain ($X = T,\,0 < \s \le 1$)
$$
\eqalign{
I_2(T) &\ll T^{-1}X^{2-2\s}\log^2X + TX^{-2\s}\log^4X + \cr&
\quad+ T^2\int_X^\infty(S^2(x) + x^{-2})x^{-1-2\s}\d x 
\cr&\ll T^{2-2\s}\log\log T.
\cr}\eqno(2.7)
$$
Combining (2.3) and (2.7), replacing $T$ by $T2^{-j}$ 
and summing all the results we finally deduce

\bigskip
THEOREM 1. {\it For $\s$ fixed we have}
$$
\int_1^T|G(\s+it)|^2\d t \ll \Bigg\{\aligned T\qquad
&(\hf < \s \le 1),\\ T\log^2T\qquad&(\s = \hf),\\
\;\;T^{2-2\s}\log T \qquad&(0 < \s < \hf).\endaligned\eqno(2.8)
$$
\bigskip\no
The lower limit of integration in (2.8) is 1 and not 0
to avoid the pole of $G(s)$ at $s=1$.
It is not difficult to see that, by using (1.5),  the validity of the last 
bound in (2.8) can be extended to the range $-1 < \s < \hf$, and 
the first bound in (2.8) to $\s > 1$
as well.  A natural problem is to try to show that for $\s =\hf$ the
integral in (2.8) is asymptotic to $CT\log^2T$.

\bigskip
\heading{\bf3. A multiple sum over zeta-zeros}
\endheading
\noindent
For a fixed $n\in\NN$, let $\g^{(1)},\ldots,\g^{(n)}$ denote ordinates
of zeta-zeros. By absolute convergence and the classical integral
$$
e^{-z} = {1\over2\pi i}\int_{(c)}w^{-z}\G(w)\d w\qquad(\R z > 0, c > 0),
$$
we have
$$
\eqalign{
\sum_{\g^{(1)}>0,\ldots,\g^{(n)}>0}e^{-\g^{(1)}\ldots\g^{(n)}/X}
&= {1\over2\pi i}\int_{(2)}
\sum_{\g^{(1)}>0,\ldots,\g^{(n)}>0}
(\g^{(1)}\ldots\g^{(n)}/X)^{-s}\G(s)\d s\cr
&= {1\over2\pi i}\int_{(2)}\G(s)G^n(s)X^s\d s.\cr}\eqno(3.1)
$$
Since $G(s)$ has a double pole at $s=1$, the function $G^n(s)$ will have
a pole of order $2n$ at $s=1$, but otherwise it is regular for $\s > -1$
and $G^n(s) \ll (1+|t|)^{4n}$ in this region. Hence by the residue theorem
and Stirling's formula for the gamma-function we obtain
$$
\eqalign{&
{1\over2\pi i}\int_{(2)}\G(s)G^n(s)X^s\d s =
X(A_{2n-1,n}\log^{2n-1}X + \cdots + A_{1,n}\log X + A_{0,n})\cr&
+ G^n(0) + {1\over2\pi i}\int_{(\e-1)}\G(s)G^n(s)X^s\d s\cr&
= X(A_{2n-1,n}\log^{2n-1}X + \cdots + A_{1,n}\log X + A_{0,n}) +
G^n(0) + O_\e(X^{\e-1}),\cr}
$$
where  $A_{2n-1,n}\not = 0,\ldots,A_{0,n}$ are effectively computable
constants. Thus we have

\bigskip
THEOREM 2. {\it For fixed $n\in\NN$ there exist effectively computable
constants $A_{2n-1,n}\not = 0,\ldots,A_{0,n}$  such that}
$$
\eqalign{
\sum_{\g^{(1)}>0,\ldots,\g^{(n)}>0}e^{-\g^{(1)}\ldots\g^{(n)}/X}&=
X(A_{2n-1,n}\log^{2n-1}X + \cdots + A_{1,n}\log X + A_{0,n})\cr&
+ G^n(0) + O_\e(X^{\e-1}),\cr}\eqno(3.2)
$$
{\it where $\g^{(1)},\ldots,\g^{(n)}$ denote ordinates
of complex zeros of $\z(s)$.}

\bigskip\noindent
If the Riemann Hypothesis holds, then the asymptotic 
formula (3.2) can be considerably sharpened. Namely
we have (see [13, eq. (14.13.8)])
$$
S_n(t) \;=\; O\left({\log t\over(\log\log t)^{n+1}}\right),
$$
where
$$
S_n(t) \;:=\;  \int_0^t S_{n-1}(u)\d u\qquad(n \ge 1,\; S_0(t)
\equiv S(t)).
$$
On the other hand, the function $f(x)$ in (1.3) admits (unconditionally)
an asymptotic expansion in terms of negative odd powers of $x$, 
in view of Stirling's formula for the gamma-function.
Thus from (1.5) we obtain, by successive integrations by parts and the above
bound for $S_n(T)$, that on the RH  the function
$G(s)$ admits analytic continuation
to $\CC$, and is of polynomial growth in $\I|s|$, provided
that $s$ stays away from its poles: the double pole at $s=1$
and simple poles at $s = -1, -3,\,\ldots\,$. As mentioned
in Section 1, these facts have been established by a different
method in Delsarte [5, p. 431].  The converse problem seems to be
interesting, namely what can be deduced about the location of zeros
of $\z(s)$ from the fact that $G(s)$ has analytic continuation to, say,
$\s > -A\,(1 < A < \infty)$?

\smallskip
It transpires that if in the above
proof we shift the line of integration (assuming RH) to $\R s = -A$,
where $A = k + \hf \ge {3\over2}$ is half of an odd natural number,
then we shall obtain in (3.2) additional main terms coming from  the poles
at $s = -1,-2, \ldots, -k$ of the integrand, 
plus an error term which will be $\ll X^{-A}$.

\bigskip\noindent
We can obtain an unconditional result  analogous to (3.2), namely 
$$
\eqalign{
\sum_{\rho^{(1)},\ldots,\rho^{(n)}}e^{-|\rho^{(1)}\ldots\rho^{(n)}|/X}&=
X(\a_{2n-1,n}\log^{2n-1}X + \cdots + \a_{1,n}\log X + a_{0,n})\cr&
+ R^n(0) + O_\e(X^{\e-1}),\cr}\eqno(3.3)
$$
where $\a_{2n-1,n}\not = 0,\ldots,\a_{0,n}$ are effectively computable
constants, $\rho^{(1)},\ldots,\rho^{(n)}$ denote complex zeros of
$\z(s)$ and, for $\s > 1$,
$$
R(s) \= \sum_\rho|\rho|^{-s} = 2\sum_{\g>0}|\rho|^{-s},
$$
and otherwise $R(s)$ is defined by analytic continuation.
This can be obtained in the region $\s > -1$ by writing
$$
\eqalign{
R(s) &= 2G(s) + 2\sum_{\g>0}(|\rho|^{-s} - \g^{-s})\cr&
= 2G(s) - 2s\sum_{\g>0}\, \int_\g^{|\rho|} x^{-s-1}\d x.\cr}\eqno(3.4)
$$
But with $\rho = \b+i\g,\,\g > 0$ we have
$$
\left|\int_\g^{|\rho|} x^{-s-1}\d x\right|
\le (|\rho|-\g)\g^{-\s-1} = (\sqrt{\b^2+\g^2} - \g)\g^{-\s-1}
\le \hf \g^{-\s-2}
$$
since $0 < \b < 1$. Hence
$$
H(s) := 2s\sum_{\g>0}(|\rho|^{-s} - \g^{-s})
$$
is regular for $\s > -1$ and in that region it satisfies
$$
H(s) \;\ll\; |s|.
$$
Therefore (3.4) provides analytic continuation of
$R(s)$ to $\s > -1$. By using the method of proof of Theorem 1
we obtain
$$
R(s) \;\ll_\e\; |t|^{1-\s+\e}\qquad(-1 < \s \le 1,\,|t| \ge t_0 > 0)
\eqno(3.5)
$$
and also
$$
\int_1^T|R(\s+it)|^2\d t \ll_{\e} \Bigg\{\aligned T^{2-2\s+\e}\qquad
&(-1<\s\le\hf),\\ T^{1+\e}\qquad&(\s \ge \hf).\endaligned\eqno(3.6)
$$
Uisng then (3.5) (or (3.6)) one obtains (3.3) similarly to the way (3.2)
was obtained.

\bigskip
\heading{\bf4. Some integrals involving $S(T)$}
\endheading
\noindent
Certain types of integrals involving the function $S(T)$ (see (1.2))
are closely related to sums over zeta-zeros,
and thus to $G(s)$. In this section we shall
investigate the evaluation of some such integrals, which do not appear
to have been treated in the literature before. We start by proving

\bigskip
THEOREM 3. {\it Let $f(t) \in C[1,T]$ satisfy}
$$
\int_1^T f^2(t)\d t \ll T\log^CT\qquad(C \ge 0).\eqno(4.1)
$$
{\it Then for fixed $r\in\NN$ we have}
$$\eqalign{&
\int_1^T S^r(t)f(t)\d t \cr&
\ll_\e \min\left(T\log^{C\over2}T(\log\log T)^{r\over2},\,
T + (\log\log T)^{{3\over2}r+\e}\int_1^T|f(t)|\d t\right).\cr}\eqno(4.2)
$$

\bigskip
{\bf Proof.} The first bound in (4.2) follows from (4.1),
the Cauchy-Schwarz inequality and the bound of K.-M. Tsang [14]
$$
\int_T^{2T}S^{2k}(t)\d t \ll T(ck)^{2k}(\log\log T)^k\quad(C > 0),
\eqno(4.3)
$$
which is uniform in $k\in\NN$. To obtain the second bound in (4.2) let,
for a given constant $\delta > 0$,
$$
H_\delta(T) := \left\{\,t\,:\,T\le t \le 2T,\,|S(t)| \ge
(\log\log T)^{{1\over2}+\delta} \right\}.
$$
Then (4.3) gives ($\mu(\cdot)$ denotes measure)
$$
\mu(H_\delta(T))(\log\log T)^{k+2k\delta} \ll T(ck)^{2k}(\log\log T)^k,
$$
and consequently
$$
\mu(H_\delta(T)) \;\ll\;T\left({ck\over(\log\log T)^\delta}\right)^{2k}.
\eqno(4.4)
$$
Choose
$$
k = \left[{1\over2c}(\log\log T)^\delta\right].
$$
Then for $T$ large enough $k \in\NN$, and (4.4) implies
$$
\mu(H_\delta(T)) \ll T2^{-2k} \le Te^{-A(\log\log T)^\delta}
\quad\left(A = {\log4\over4c}\right).\eqno(4.5)
$$
Thus if $\delta >1$, then for any fixed $C_1 > 0$ we have from (4.5)
$$
\mu(H_\delta(T)) \ll T(\log T)^{-C_1}.\eqno(4.6)
$$
Now suppose that $\delta > 1$. Then using (1.2) and (4.6) we have
$$
\eqalign{
&\int_T^{2T}S^r(t)f(t)\d t = \int_{H_\delta(T)} +
\int_{[T,2T]\setminus{H_\delta(T)}}\cr&
\ll \left(\int_{H_\delta(T)} S^{2r}(t)\d t\right)^{1/2}
\left(\int_T^{2T}f^2(t)\d t
\right)^{1/2}\cr&
+ (\log\log T)^{r({1\over2}+\delta)}\int_T^{2T}|f(t)|\d t\cr&
\ll (T(\log T)^{2r-C_1})^{1/2}(T\log^CT)^{1/2} + 
(\log\log T)^{r({1\over2}+\delta)}\int_T^{2T}|f(t)|\d t\cr&
\ll T + (\log\log T)^{{3r\over2}+\e}\int_T^{2T}|f(t)|\d t
\cr}
$$
with $C_1 = 2r + C,\,\delta = 1 + \e/r$. Replacing $T$ by $T2^{-j}\,
(j \in \NN)$ and adding up the resulting estimates we complete the
proof of (4.2).

\medskip
The integrals which seem of interest are e.g.,
$$
\int_1^T S(t)|\zt|^2\d t,\quad  \int_1^T S^2(t)|\zt|^2\d t\eqno(4.7)
$$
and 
$$
\int_1^T |\zt|^2\d S(t). \eqno(4.8)
$$
An integration by parts shows that the integral in (4.8) equals
$$
|\zt|^2S(t)\Big|_1^T - 2\int_1^TS(t)Z(t)Z'(t)\d t,
$$
where Hardy's function $Z(t)$ (see [8], [11]) is a real-valued function of $t$
satisfying $|Z(t)| = |\zt|$, and given by
$$
Z(t) := \zt\chi^{-1/2}(\hf+it),\quad \chi(s) = {\z(s)\over\zeta(1-s)}
= 2^s\pi^{s-1}\sin(\hf \pi s)\G(1-s).
$$
Since
$$
\int_0^T|Z(t)Z'(t)|\d t \le \left(\int_0^T|\zt|^2\d t\right)^{1/2}
\left(\int_0^T|Z'(t)|^2\d t\right)^{1/2} \ll T\log^2T,
$$
it follows on using Theorem 3 that
$$
\int_1^T |\zt|^2\d S(t) \ll_\e T\log^2T(\log\log T)^{{3\over2}+\e}.\eqno(4.9)
$$
Similarly we have
$$
\int_1^T |\zt|^2S(t)\d t \ll_\e T\log T(\log\log T)^{{3\over2}+\e},
\eqno(4.10)
$$
and
$$
\int_1^T |\zt|^2S^2(t)\d t \ll_\e T\log T(\log\log T)^{3+\e}.\eqno(4.11)
$$
The bounds (4.9)--(4.11) appear to be, at present, the strongest
unconditional bounds that can be obtained.

On the other hand, the above integrals can be related to sums over
zeta-zeros. For example,  the integral in (4.8) is
$$
\eqalign{
\int_1^T |\zt|^2\d N(t) &- \int_1^T{1\over2\pi}\log{t\over2\pi}\cdot
|\zt|^2\d t + O(1) \cr&
= \sum_{0<\g\le T}|\z(\hf+i\g)|^2 - {T\over2\pi}\log^2T + O(T\log T).
\cr}
$$
This gives, on using (4.9),
$$\eqalign{
\sum_{0<\g\le T}|\z(\hf+i\g)|^2 &= {T\over2\pi}\log^2T + O(T\log T)
+ \int_1^T |\zt|^2\d S(t)\cr&
\ll_\e  T\log^2T(\log\log T)^{{3\over2}+\e}.\cr}\eqno(4.12)
$$
We recall the standard notation (see [8] and [9])
$$
\int_0^T|\zt|^2\d t = T\log\left({T\over2\pi}\right) + (2C_0 - 1)T + E(T),
$$
where $C_0$ denotes Euler's constant. Then by using integration
by parts, (1.6) and the bound $E(T) \ll T^c$ with suitable
$c < 1/3$ (see [8]) we have
$$
\eqalign{&
\int_0^T S(t)|\zt|^2\d t = \int_0^TS(t)\left(\log{t\over2\pi}+2C_0 + E'(t)
\right)\d t\cr&
= O(\log^2T) + \int_0^TS(t)E'(t)\d t = O(T^{1/3}) - \int_0^T E(t)\d S(t)
\cr&
= O(T^{1/3}) - \int_0^T E(t)\left(\d N(t) - {1\over2\pi}\log{t\over2\pi}
\d t + \d O\bigl({1\over t}\bigr)\right)
\cr&
= - \sum_{0<\g\le T}E(\g) + O(T^{1/3}) + {1\over2\pi}\int_0^T E(t)
\log{t\over2\pi}\d t.
\cr}
$$
The last integral equals
$$
\eqalign{
\int_0^T(E(t) - \pi + \pi)\log{t\over2\pi}\d t &= O(T^{3/4}\log T)
+ \pi \int_0^T\log{t\over2\pi}\d t\cr&
= \pi T\log T + O(T),
\cr}
$$
since we have (see [9])
$$
\int_0^T(E(t) - \pi)\d t \ll T^{3/4}.
$$
Therefore by using (4.10) we obtain
$$
\eqalign{
\sum_{0<\g\le T}E(\g) &= \hf T\log T + O(T) - \int_0^T S(t)|\zt|^2\d t\cr&
\ll_\e T\log T(\log\log T)^{{3\over2}+\e}.
\cr}\eqno(4.13)
$$
A similar calculation will also give (see [9] and [10])
$$
\sum_{0<\g\le T}E_2(\g) \ll T^{3/2}\log T,
\quad
\sum_{0<\g\le T}E_2(\g) = \Omega_\pm(T^{3/2}\log T),\eqno(4.14)
$$
where $E_2(T)$ is the error term in the asymptotic formula for
the fourth power of $|\zt|$.

\bigskip
The importance of the sum
$$
\sum_{0<\g\le T}|\z(\hf + i\g)|^2\eqno(4.15)
$$
lies in the fact that it  identically vanishes if the Riemann Hypothesis
 holds. The
unconditional bound (4.12) seems to be very weak. However this reflects
the enormous difficulty of settling the Riemann Hypothesis.
It may be remarked that
a more general sum than the one in (4.15) was treated by S.M. Gonek [7]. 
He proved, assuming the Riemann Hypothesis, that
$$
\sum_{0<\g\le T}\left|\z\left({1\over2}+i\left(\g + {\a\over L}\right)\right)
\right|^2
= \left(1 - \left({\sin\pi\a\over\pi\a}\right)^2\right){T\over2\pi}
\log^2T + O(T\log T)\eqno(4.16)
$$
holds uniformly for $|\a| \le \hf L$, where $L = {1\over2\pi}
\log({T\over2\pi})$. It would be interesting to recover this result
unconditionally, but our method of proof  does not
seem capable of achieving this. 

\medskip
One can treat the integrals in (4.9)-(4.11) by  using Lemma 2
of Bombieri-Hejhal [1], which (after taking the imaginary part)
provides an explicit expression for $S(T)$. The best this could
give (in view of $O(1)$ in the error term) for the integral in (4.9)
is the bound $O(T\log^2T)$, which is still quite weak. Assuming the RH
Lemma 2 of [1] will yield
$$
\int_0^T S(t)|\zt|^2\d t \;=\;O(T\log T).\eqno(4.17)
$$
It remains elusive whether the bound in (4.17) gives the correct
order of magnitude for the integral on the left-hand side. Is the
integral $\Omega_\pm(T\log T)$? This seems to be difficult to settle, even
if the Riemann Hypothesis is assumed.

\vskip1cm
{\bf Acknowledgement}. I wish to thank Prof. Akio Fujii for valuable remarks.
\vglue2cm
\vfill\break

\Refs

\bigskip
\item{[1]} E. Bombieri and D.A. Hejhal, {\it On the distribution of zeros
of linear combinations of Euler products}, Duke Math. J. {\bf80}(1995),
821-862.

\item{[2]} I.C. Chakravarty, {\it The secondary zeta-functions},
Journal Math. Anal. Appl. {\bf30}(1970), 280-294.

\item{[3]} I.C. Chakravarty, {\it Certain properties
of a pair of secondary zeta-functions},
Journal Math. Anal. Appl. {\bf35}(1971), 484-495.

\item{[4]} H. Davenport, {\it Multiplicative Number Theory }
2nd edition, GTM{\bf74}, Springer, New York-Heidelberg-Berlin, 1980.

\item{[5]} J. Delsarte, {\it Formules de Poisson avec reste},
Journal Anal. Math. {\bf17}(1966), 419-431.

\item{[6]} A. Fujii, {\it The zeros of the zeta function and Gibbs's 
phenomenon}, Comment. Math. Univ. Sancti Pauli {\bf32}(1983), 229-248.

\item{[7]} S.M. Gonek, {\it Mean values of the Riemann zeta-function and
its derivatives}, Invent. math. {\bf75}(1984), 123-141.

\item{[8]} A. Ivi\'c, {\it The Riemann zeta-function}, John Wiley
\& Sons, New York, 1985.

\item{[9]} A. Ivi\'c, {\it The mean values of the Riemann zeta-function}, 
Tata Institute of Fundamental Research, Lecture Notes {\bf82}, 
    Bombay 1991 (distr. Springer Verlag, Berlin etc.).

\item{[10]} A. Ivi\'c, {\it On the error term for the fourth moment 
of the Riemann zeta-function},
Journal London Math. Society {\bf60}(2)(1999), 21-32.

\item{[11]} A.A. Karatsuba and S.M. Voronin, {\it The Riemann zeta-function},
Walter de Gruyter, Berlin-New York, 1992.

\item{[12]} E.C. Titchmarsh, {\it Introduction to the Theory of Fourier
Integrals},  Clarendon Press, Oxford, 1948.

\item{[13]} E.C. Titchmarsh, {\it The theory of the Riemann
zeta-function}, 2nd edition, Oxford University Press, Oxford, 1986.

\item{[14]} K.-M. Tsang, {\it Some $\Omega$-theorems for the Riemann
zeta-function}, Acta Arith. {\bf46}(1986), 369-395.

\bigskip
\bigskip

Aleksandar Ivi\'c

Katedra Matematike RGF-a

Universitet u Beogradu

\DJ u\v sina 7, 11000 Beograd, Serbia

{\tt aivic\@rgf.bg.ac.yu, aivic\@matf.bg.ac.yu}

\endRefs

\bye